\documentclass[11pt]{article}
\usepackage{amsmath,amsfonts,amssymb,amsthm} 
\usepackage[explicit]{titlesec}
\usepackage{graphicx}
\usepackage{enumitem}
\usepackage{subcaption}
\usepackage{setspace}
\setlist{nolistsep}
\usepackage{stmaryrd}
\usepackage{mathpple}
\usepackage{enumerate}
\usepackage{textcomp}
\usepackage[utf8]{inputenc}
\usepackage{mathtools}
\usepackage{epsfig}
\usepackage{xcolor}
\usepackage{array}
\usepackage{fancyhdr}
\usepackage[colorlinks=false]{hyperref}
\usepackage[noabbrev]{cleveref}
\usepackage{dsfont}
\usepackage{wrapfig}
\usepackage{float}
\usepackage[font={small}]{caption}
\usepackage[left=3cm,right=3cm,top=2.5cm,bottom=2.5cm]{geometry}

\DeclareCaptionFont{fig}{\bfseries\footnotesize\fontfamily{phv}\selectfont}
\DeclareMathAlphabet{\mathcal}{OMS}{cmsy}{m}{n}

\newcommand{\IP}{\mathbb{P}}

\newcommand{\N}{\mathbb{N}}

\newcommand{\E}{\mathbb{E}}

\newcommand{\IL}{\mathbb{L}}

\newcommand{\cB}{\mathcal{B}}
\newcommand{\cF}{\mathcal{F}}

\newcommand{\cU}{\mathcal{U}}
\newcommand{\cR}{\mathcal{R}}
\newcommand{\cN}{\mathcal{N}}
\newcommand{\cM}{\mathcal{M}}

\newcommand{\1}{\mathds{1}}

\newcommand{\eps}{\varepsilon}
\newcommand{\cMc}{\langle M \rangle}

\newcommand{\floor}[1]{\lfloor #1 \rfloor}

\newcommand{\comb}[2]{\Big({\scriptstyle\begin{matrix}#1 \\ #2\end{matrix}}\Big)}

\font\calcal=cmsy10 scaled\magstep1
\def\build#1_#2^#3{\mathrel{\mathop{\kern 0pt#1}\limits_{#2}^{#3}}}
\def\liml{\build{\longrightarrow}_{}^{{\mbox{\calcal L}}}}

\numberwithin{equation}{section}

\setitemize[1]{label=--,partopsep=\parskip,topsep=-\parskip,itemsep=0pt,parsep=0pt}
\titleformat{\chapter}{\fontfamily{phv}\selectfont\LARGE\bfseries}{\thesection}{1em}{\fontfamily{phv}\selectfont\LARGE\bfseries #1}
\titleformat{\section}{\centering\fontfamily{phv}\selectfont\bfseries}{\thesection}{1em}{\centering\fontfamily{phv}\selectfont\bfseries #1}
\titleformat{\subsection}{\fontfamily{phv}\selectfont\small\bfseries\centering}{\thesubsection}{1em}{\fontfamily{phv}\selectfont\small\bfseries #1}
\titleformat{\subsubsection}{\fontfamily{phv}\selectfont\footnotesize\bfseries\centering}{\thesubsubsection}{1em}{\fontfamily{phv}\footnotesize\selectfont\bfseries #1}

\renewcommand*\thesection{\arabic{section}}

\makeatletter
\renewcommand{\@seccntformat}[1]{\llap{{\csname the#1\endcsname}\hspace{1em}}}   

\newtheoremstyle{thmit}{10pt}{10pt}
{\normalfont\itshape}
{}
{\small\bf\fontfamily{phv}\selectfont}
{\;}
{0.25em}
{\small\fontfamily{phv}\selectfont\thmname{#1}\nobreakspace\footnotesize\thmnumber{#2}.
\thmnote{\nobreakspace\the\thm@notefont\fontfamily{phv}\selectfont\footnotesize\bfseries\--\nobreakspace#3}}

\newtheoremstyle{rmq}
{5pt}
{5pt}
{\normalfont}
{}
{\small\bf\fontfamily{phv}\selectfont}
{\;}
{0.25em}
{\small\fontfamily{phv}\selectfont\thmname{#1}\nobreakspace\footnotesize\thmnumber{\@ifnotempty{#1}{}\@upn{#2}}
\thmnote{\nobreakspace\the\thm@notefont\fontfamily{phv}\selectfont\footnotesize\bfseries--\nobreakspace#3.}}
\makeatother

\newcounter{count}
\numberwithin{count}{section}
\newcounter{alpha}

\theoremstyle{thmit}
\newtheorem{theorem}[count]{Theorem\small}

\newtheorem{lemma}[count]{Lemma\small}

\newtheorem{remark}[count]{Remark\small}

\theoremstyle{rmq}

\renewenvironment{proof}[2]
{\paragraph{\fontfamily{phv}\selectfont\bfseries\small Proof of #1 \footnotesize #2.}}%
{\begin{flushright}
\qed
\end{flushright}}

\begin{document}
\pagestyle{fancy}

\title{\vspace{-2ex}
\fontfamily{phv}\selectfont\bfseries\Large
Introducing smooth amnesia to the memory of the Elephant Random Walk}
\author{\fontfamily{phv}\selectfont\bfseries
Lucile Laulin}
\date{}
\AtEndDocument{\bigskip{\footnotesize%
\textsc{Université de Bordeaux, Institut de Mathématiques de Bordeaux,
UMR 5251, 351 Cours de la Libération, 33405 Talence cedex, France.} \par  
\textit{E-mail adress :} \href{mailto:lucile.laulin@math.u-bordeaux.fr}{\texttt{lucile.laulin@math.u-bordeaux.fr}} \par
}}

\maketitle

\begin{center}
\begin{minipage}[c]{0.8\textwidth}
{\small\section*{Abstract}\vspace{-1ex}
 This paper is devoted to the asymptotic analysis of the amnesic elephant random walk (AERW) using a martingale approach. More precisely, our analysis relies on asymptotic results for multidimensional martingales with matrix normalization.
 In the diffusive and critical regimes, we establish the almost sure convergence and the quadratic strong law for the position of the AERW. The law of iterated logarithm is given in the critical regime. The distributional convergences of the AERW to Gaussian processes are also provided. In the superdiffusive regime, we prove the distributional convergence as well as the mean square convergence of the AERW. 
}\medskip\\
\small
{\bf MSC:} primary 60G50; secondary 60G42; 60F17\medskip\\
{\bf Keywords : }Elephant random walk; Amnesic random walk; Multi-dimensional martingales; Almost sure convergence; Asymptotic normality; Distributional convergence\end{minipage}   
\end{center}
\setlength{\parindent}{0pt}

\section{Introduction}
The Elephant Random Walk (ERW) is a discrete-time random walk, introduced by Schütz and Trimper \cite{Schutz2004} in the early 2000s. At first, the ERW was used in order to see how long-range memory affects the random walk and induces a crossover from a diffusive to superdiffusive behavior. It was referred to as the ERW in allusion to the traditional saying that elephants can always remember anywhere they have been.  
The elephant starts at the origin at time zero, $S_0 = 0$. At time $n = 1$,  the elephant moves one step to the right with probability $q$ and to the left with probability $1-q$ for some $q$ in $[0,1]$.  
Afterwards, at time $n+1$, the elephant chooses uniformly at random an integer $k$ among the previous times $1,\ldots,n$. Then, it moves exactly in 
the same direction as that of time $k$ with probability $p$ or the opposite direction with the probability $1-p$, where the parameter $p$ 
stands for the memory parameter of the ERW.
The position of the elephant at time $n+1$ is given by
\begin{equation} 
\label{POS-bERW}
S_{n+1} = S_n + X_{n+1}
\end{equation}
where $X_{n+1}$ is the $(n+1)$-th increment of the random walk, such that
\begin{equation}
\label{ERW-X-a-b}
X_{n+1}=\alpha_{n+1} X_{\beta_{n+1}}
\end{equation}
where $\alpha_{n+1}\sim \cR(p)$ and $\beta_{n+1}\sim\cU(1,n)$ are mutually independent and independant of the past.
The ERW shows three differents regimes depending on the location of its memory parameter $p$ with respect to the critical value $p=3/4$.

On the one hand, a wide literature is now available on the ERW in dimension $d=1$ thanks to a variety of approaches. Baur and Bertoin \cite{Baur2016} used the connection to P\'{o}lya-type urns as well as functional limit theorems for multitype branching processes due to Janson \cite{Janson2004}. Bercu \cite{Bercu2018} and Coletti et al. \cite{Coletti2017} used martingales to obtain the almost sure convergence and asymptotic normality, among other results. Kürsten \cite{Kursten2016} and Businger \cite{Businger2018} used the construction of random trees with Bernoulli percolation. A strong law of large numbers and a central limit theorem for the position of the ERW, properly normalized, were established in the diffusive regime $p< 3/4$ and the critical regime $p=3/4$, see \cite{Baur2016,Bercu2018,Coletti2017,Coletti2019,Vazquez2019}. In the superdiffusive regime $p>3/4$, Bercu \cite{Bercu2018} proved that the limit of the position of the ERW is not Gaussian and Kubota and Takei \cite{Kubota2019} showed that the fluctuation of the ERW around this limit is Gaussian. 

On the other hand, over the last years, various processes derivated from the ERW have recevied a lot of attention. Bercu and Laulin in \cite{BercuLaulin2019} extended all the results of \cite{Bercu2018} to the multi-dimensional ERW (MERW) where $d \geq 1$ and to its center of mass \cite{LaulinBercu2020} using a martingale approach, while Bertenghi used the connection \cite{Bertenghi2020} to Pólya-type urns for the MERW. The ERW with stops or minimal RW, changing in particular the distribution of $\alpha_{n}$, has also been investigated \cite{Bercu2021,Bercu2022,Gut2019,Takei2020}. The ERW with reinforced memory has been studied by Baur \cite{baur2019} via the urn approach, and Laulin \cite{Laulin2021} using martingales.
\medskip

The idea of this paper is to use the approach developped in \cite{LaulinBercu2020} and \cite{Laulin2021} to study how changing the memory allows us to induce amnesia to the ERW.
More precisely, the distribution of the memory $\beta_n$ of our new variation of the ERW is such that the probability of choosing a fixed instant $k\in\N^*$ at time $n\geq k$ decreases approximatly with speed $(k/n)^{\beta}$ for some amnesia parameter $\beta\geq 0$.
\smallskip

The very interesting question of amnesic elephant random walk (AERW) has not been investigated a lot. Gut and Stadmüller \cite{GutA2021,GutStad2021} studied variations of the memory for the special cases of ERW with delays or gradually increasing memory. In \cite{GutA2021} the elephant could stop and only remember the first (and second) step it tooks. Consequently, it did not induced a phase transition. In \cite{GutStad2021}, the elephant only remembered a portion of its past (recent or distant), this portion being fixed or depending on the time $n$, but was always ``small".
\smallskip

The entire study we conduct below can be generalized when $\beta<0$ is not an integer. This can be interpreted as cases where the elephant remembers more vividly the first steps it performed. When $\beta<-1$, it appears that the AERW only have one regime that is the diffusive regime. This observation is coherent with the work of Gut and Stadmüller \cite{GutStad2021}.
\smallskip 

\begin{figure}
\begin{subfigure}{0.3\linewidth}
    \includegraphics[width=5cm]{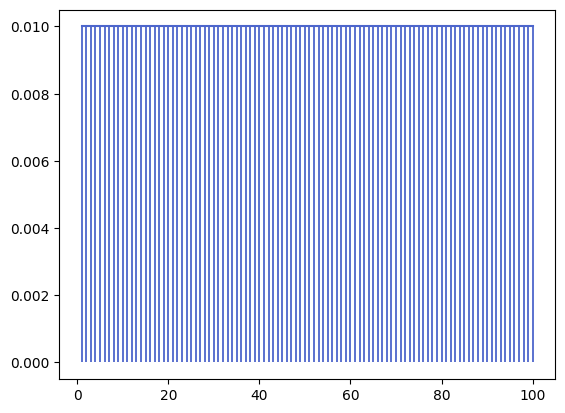}
    \caption{$\beta =0$}
    \label{AERW-picture-beta0}
\end{subfigure}\hfill
\begin{subfigure}{0.3\linewidth}
    \includegraphics[width=5cm]{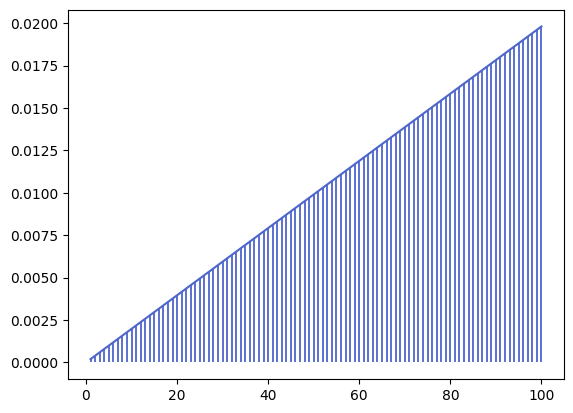}
    \captionof{figure}{$\beta =1$}
    \label{AERW-picture-beta1}
\end{subfigure}
\hfill
\begin{subfigure}{0.3\linewidth}
    \includegraphics[width=5cm]{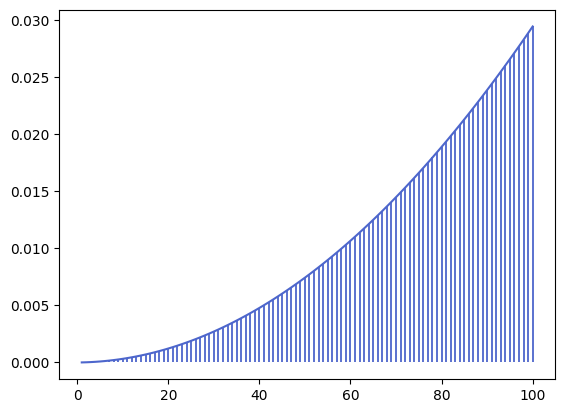}
    \captionof{figure}{$\beta =2$}
    \label{AERW-picture-beta2}
\end{subfigure}\\
\begin{subfigure}{0.3\linewidth}
    \includegraphics[width=5cm]{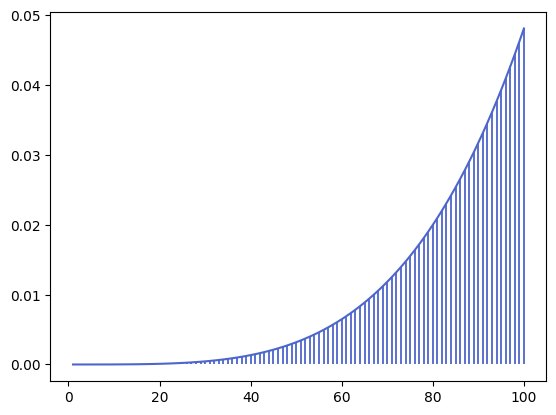}
    \captionof{figure}{$\beta =3$}
    \label{AERW-picture-beta3}
\end{subfigure}\hfill
\begin{subfigure}{0.3\linewidth}
    \includegraphics[width=5cm]{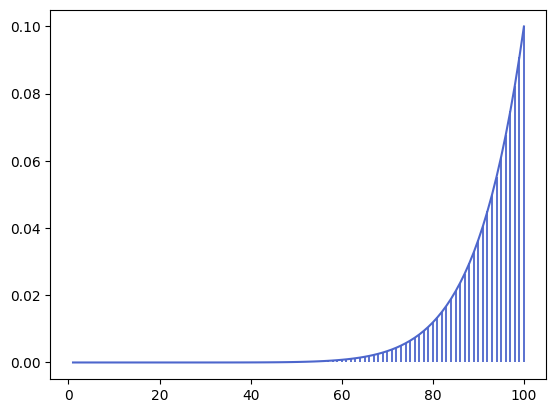}
    \captionof{figure}{$\beta =10$}
    \label{AERW-picture-beta10}
\end{subfigure}
\hfill
\begin{subfigure}{0.3\linewidth}
    \includegraphics[width=5cm]{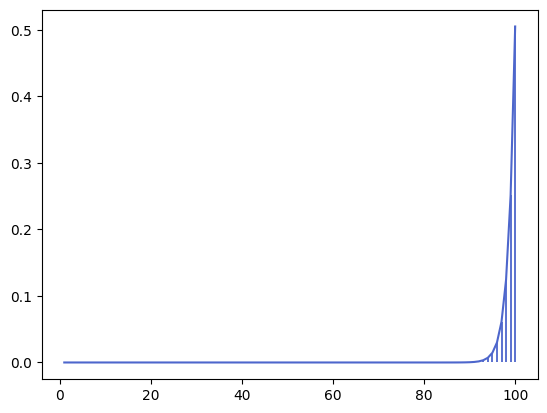}
    \captionof{figure}{$\beta =100$}
    \label{AERW-picture-beta100}
\end{subfigure}
\caption{Mass function of the memory depending on the value of $\beta$.}
\label{fig-beta}
\end{figure}

The AERW will appear to be non-Markovian, as the reinfroced ERW. However, unlike the reinforced ERW, the AERW can not be studied using Pólya-type urns.
The major change for the AERW is that the distribution of the memory $\beta_n$ in equation \eqref{ERW-X-a-b} is no longer uniform but depends on the amnesia parameter $\beta\geq 0$. In this approach, the elephant chooses an instant according to $\beta_{n+1}$ as follows,
\begin{equation}
\label{beta-distribution}
    \IP(\beta_{n+1}=k)= \frac{(\beta+1)\Gamma(k+\beta)\Gamma(n)}{\Gamma(k)\Gamma(n+\beta+1)} = \frac{(\beta+1)}{n}\frac{\mu_k}{\mu_{n+1}}
    \quad
    \text{for $1\leq k\leq n$},
\end{equation}
where 
\begin{equation}
    \label{AERW-DEF-muN}
    \mu_n = \prod_{k=1}^{n-1}\Big(1+\frac{\beta}{k}\Big)= \frac{\Gamma(n+\beta)}{\Gamma(n)\Gamma(\beta+1)}.
\end{equation}
The case $\beta=0$ corresponds to the traditionnal ERW. As $\beta$ grows, the probability of choosing a recent instant gets bigger, see the illustrative Figure \ref{fig-beta}.

We have by definition of the step $X_{n+1}$ given in \eqref{ERW-X-a-b} and the distribution $\beta_{n+1}$ \eqref{beta-distribution} that
\begin{align}
\label{AERW-ESPC-X}
    \E[X_{n+1}\mid\cF_n]&=\E[\alpha_{n+1}]\E[X_{\beta_{n+1}}\mid\cF_n]\nonumber \\
    &=(2p-1)\E\Big[\sum_{k=1}^n X_k\1_{\beta_{n+1}=k}\mid\cF_n\Big]\nonumber\\
    &=\frac{(2p-1)(\beta+1)}{n\mu_{n+1}}\sum_{k=1}^n X_k  \mu_{k}.
\end{align}
Then, denote $a=2p-1$ and
\begin{equation}
    Y_n = \sum_{k=1}^n X_k  \mu_k.
\end{equation}
We deduce from \eqref{AERW-ESPC-X} that
\begin{equation}
\label{AERW-ESPC-Y}
    \E[Y_{n+1}\mid \cF_n] = \Big(1 + \frac{a(\beta+1)}{n}\Big)Y_n.
\end{equation}
Hereafter, for any $n\geq 1$, let
\begin{equation}
\label{AERW-DEF-nuN}
    \gamma_n = 1+\frac{a(\beta+1)}{n}
\end{equation}
and 
\begin{equation}
\label{AERW-DEF-alphaN}
a_{n}=\prod_{k=1}^{n-1}\gamma_k^{-1} =\frac{\Gamma(n)\Gamma(a(\beta+1)+1)}{\Gamma(n+a(\beta+1))}.
\end{equation}
It follows from standard resultats on the Gamma function that
\begin{equation}
\label{AERW-lim-alphaN}
    \lim_{n\to\infty} n^{a(\beta+1)}a_n =  \Gamma(a(\beta+1)+1)
\end{equation}
and
\begin{equation}
\label{AERW-lim-muN}
    \lim_{n\to\infty} n^{-\beta}\mu_n = \Gamma(\beta+1).
\end{equation}
Our strategy for proving asymptotic results for the AERW is as follows. On the one hand, the behavior of the position $S_n$ is closely
related to the one of the sequences $(M_n)$ and $(N_n)$ defined, for all $n\geq 0$, by
\begin{equation}
\label{AERW-DEF-M-N}
    M_n =  a_n Y_n \quad\text{and}\quad N_n = S_n + \frac{a(\beta+1)}{\beta-a(\beta+1)}\mu_n^{-1}Y_n.
\end{equation}
We immediatly get from \eqref{AERW-ESPC-Y} and \eqref{AERW-DEF-alphaN} that $(M_n)$ is a locally square-integrable martingale adapted to $(\cF_n)$. Moreover, we have from \eqref{AERW-ESPC-X} that
\begin{equation*}
    E\Big[S_{n+1}+\frac{a(\beta+1)}{\beta-a(\beta+1)}\mu_{n+1}^{-1}Y_{n+1}\mid \cF_n\Big]=S_n + \frac{a(\beta+1)}{\beta-a(\beta+1)}\mu_n^{-1}Y_n
\end{equation*}
which also means that $(N_n)$ is also a locally square-integrable martingale adapted to $\cF_n$.
On the other hand, we can rewrite $S_n$ as 
\begin{equation}
\label{AERW-SN-MG}
    S_n = N_n - \frac{a(\beta+1)}{\beta-a(\beta+1)}(\mu_n a_n)^{-1}M_n
\end{equation}
and equation \eqref{AERW-SN-MG} allows us to establish the asymptotic behavior of the AERW via an extensive use of the martingale theory.

\smallskip
The main results of this paper are given in Section \ref{AERW-S-main-results}. We first investigate the diffusive regime and we establish the strong law of large numbers, the law of iterated logarithm and the quadratic strong law for the AERW. The functional central limit theorem is also provided. Next, we prove similar results in the critical regime. Finally, we establish a strong limit theorem in the superdiffusive regime. Our martingale approach is described in Section \ref{AERW-S-martingale-approach}. Finally, we give some of the technical proofs in Section \ref{AERW-S-proofs}.

\section{Main results}
\label{AERW-S-main-results}
\subsection{The diffusive regime}
Our first result deals with the strong law of large numbers for the AERW in the diffusive regime where $p<\frac{4\beta+3}{4(\beta+1)}$.
The following strong law for the AERW will be deduced from both the strong laws for $(N_n)$ and $(M_n)$.
\begin{theorem}
\label{AERW-THM-CVPS-DIF}
We have the almost sure convergence
\begin{equation}
    \lim_{n\to\infty} \frac{S_n}{n}=0 \quad\text{a.s.}
\end{equation}
\end{theorem}
The almost sure rate of convergence for the AERW is as follows, for
\begin{equation*}
    \sigma_\beta^2 = \frac{2\beta+1-a}{(1-a)(1+2\beta-2a(\beta+1))}.
\end{equation*}
\begin{theorem}
\label{AERW-THM-LFQ-DIF}
We have the quadratic strong law
\begin{equation}
\label{AERW-QUAD-DIF}
    \lim_{n \rightarrow \infty} \frac{1}{\log n} \sum_{k=1}^n \frac{S_k^2}{k^2}
        =\sigma_\beta^2\quad  \text{a.s.}
\end{equation}
\end{theorem}
Hereafter, we are interested in the distributional convergence of the AERW, which holds in the Skorokhod space $D([0,\infty[)$ of right-continuous functions with left-hand limits.
\begin{theorem}
\label{AERW-CLFT-AERW-DR}
The following convergence in distribution in $D([0,\infty[)$ holds
\begin{equation}
\label{AERW-CVL-cont-DIF}
\Big(\frac{S_{\floor{nt}}}{\sqrt{n}}, \ t\geq 0\Big) \Longrightarrow \big(W_t, \ t\geq 0\big)
\end{equation}
where $\big(W_t, \ t\geq 0\big)$ is a real-valued centered Gaussian process starting from the origin with covariance
\begin{align}
\E[W_sW_t] & = \frac{a(1+\beta)(1 - a)+a\beta}{(2(\beta+1)(1-a)-1)(a-\beta(1-a))(1-a)}s\Big(\frac{t}{s}\Big)^{a-\beta(1-a)} \nonumber \\
& \quad+ \frac{\beta}{(\beta(1-a)-a)(1-a)}s
\end{align}
for $0<s\leq t$. In particular, we have
\begin{equation}
\label{AERW-CVL-DIF}
\frac{S_n}{\sqrt{n}} \liml \cN \left(0, \sigma_\beta^2\right).
\end{equation}
\end{theorem}
\begin{remark}
When $\beta=0$ we find again the results from \cite{Baur2016} for the ERW
\begin{equation*}
\Big(\frac{S_{\floor{nt}}}{\sqrt{n}}, \ t\geq 0\Big) \Longrightarrow \big(W_t, \ t\geq 0\big)
\end{equation*}
where $\big(W_t, \ t\geq 0\big)$ is a real-valued mean-zero Gaussian process starting from the origin and
\begin{equation*}
\label{AERW-ESP-W}
\E[W_sW_t] = \frac{1}{1-2a}s\Big(\frac{t}{s}\Big)^{a}.
\end{equation*}
\end{remark}
\subsection{The critical regime}
Hereafter, we investigate the critical regime where $p=\frac{4\beta+3}{4(\beta+1)}$. It is interesting to notice that, when $\beta$ is really large (or $\beta \to \infty$) the critical regime is reached for the memory parameter $p=1$. Hence, the greater $\beta$ is, the more there are values of the memory parameter $p$ for which the AERW stays in the diffusive regime; but whatever the value of $\beta$, we still observe a phase transition.

\begin{theorem}
\label{AERW-THM-CVPS-CR}
We have the almost sure convergence
\begin{equation}
    \lim_{n\to\infty} \frac{S_n}{\sqrt{n}\log n}=0\quad\text{a.s.}
\end{equation}
\end{theorem}
The almost sure rates of convergence for the AERW are as follows.
\begin{theorem}
\label{AERW-THM-LFQ-LIL-CR}
We have the quadratic strong law
\begin{equation}
\label{AERW-QUAD-CR}
    \lim_{n \rightarrow \infty} \frac{1}{\log \log n} \sum_{k=1}^n \frac{S_k^2}{(k\log k)^2}=(2\beta+1)^2\quad  \text{a.s.}
\end{equation}
In addition, we also have the law of iterated logarithm
\begin{equation}
\label{AERW-LIL-CR}
 \limsup_{n \rightarrow \infty} \frac{S_n^2}{2 n \log n \log \log \log n}= (2\beta+1)^2  \quad \text{a.s.}
\end{equation}
\end{theorem}
Once again, our next result concerns the asymptotic normality of the AERW.
\begin{theorem}
\label{AERW-THM-CLFT-CR}
The following convergence in distribution in $D([0,\infty[)$ holds
\begin{equation}
\label{AERW-CLFT-CR}
\Big(\frac{S_{\floor{n^t}}}{\sqrt{n^t\log n}}, \ t\geq 0\Big) \Longrightarrow (2\beta+1)\big(B_t, \ t\geq 0\big)
\end{equation}
where $(B_t,\ t\geq0)$ is a one-dimensional standard Brownian motion.
In particular, we have the asymptotic normality
\begin{equation}
\frac{S_n}{\sqrt{n\log n}} \liml \cN \Big(0, (2\beta+1)^2 \Big).
\end{equation}
\end{theorem}
\subsection{The superdiffusive regime}
Finally, we focus our attention on the superdiffusive regime where $p>\frac{4\beta+3}{4(\beta+1)}$.
\begin{theorem}
\label{AERW-THM-CVPS-SDIF}
We have the following distributional convergence in $D([0,\infty[)$
\begin{equation}
\label{AERW-CVPS-cont-SDIF}
\Big(\frac{S_{\floor{nt}}}{n^{a(\beta+1)}}, \ t\geq 0\Big) \Longrightarrow (\Lambda_t, \ t\geq 0)
\end{equation}
where the limiting $\Lambda_t = t^{a(\beta+1)}L_\beta$, $L_\beta$ being some non-denegerate random variable. In particular, we have
\begin{equation}
\label{AERW-CVPS-SDIF}
    \lim_{n\to\infty} \frac{S_n}{n^{a(\beta+1)-\beta}}=L_\beta \quad\text{a.s.}
\end{equation}
where the limiting $L_\beta$ is a non-degenerate random variable. We also have the mean square convergence
\begin{equation}
\label{AERW-CVL2-SDIF}
    \lim_{n\to\infty} \E\Big[\Big|\frac{S_n}{n^{a(\beta+1)-\beta}}-L_\beta\Big|^2\Big]=0.
\end{equation}
\end{theorem}
\begin{remark}
\label{AERW-THM-L2-SDIF}
The expected value of $L_\beta$ is 
\begin{equation}
\label{AERW-ESP-SDIF}
    \E[L_\beta]=  \frac{a(\beta+1)(2q-1)\Gamma(\beta+1)}{\big(a(\beta+1)-\beta\big)\Gamma\big(a(\beta+1)+1\big)}
\end{equation} 
while its second order moment is given by
\begin{equation}
\label{AERW-L2-SDIF}
    \E\big[L_\beta^2\big]= \frac{a^2(\beta+1)^2\Gamma(\beta+1)^2\Gamma\big(2(a-1)(\beta+1)+1\big)}{\big(a(\beta+1)-\beta\big)^2\Gamma\big((2a-1)(\beta+1)+1\big)^2}.
\end{equation}
    When $\beta=0$ we find again the expected values for the ERW from \cite{Bercu2018}
    \begin{equation*}
        \E[L] = \frac{2q-1}{\Gamma(a+1)}\quad\text{and}\quad
        \E[L^2] = \frac{1}{(2a-1)\Gamma(2a)}.
    \end{equation*}
\end{remark}
\section{A two-dimensional martingale approach}
\label{AERW-S-martingale-approach}
In order to investigate the asymptotic behavior of $(S_n)$, we introduce the two-dimensional martingale $(\cM_n)$ defined by
\begin{equation}
\label{AERW-DEF-MMN}
    \cM_n = \begin{pmatrix} N_n \\ M_n\end{pmatrix}
\end{equation}
where $(M_n)$ and $(N_n)$ are the two locally square-integrable martingales introduced in \eqref{AERW-DEF-M-N}. As for the CMERW and the RERW, the main difficulty we face is that the predictable quadratic variations of $(M_n)$ and $(N_n)$ increase to infinity with two different speeds. A matrix normalization will again be necessary to establish the asymptotic behavior of the AERW. We will alternatively study $(\cM_n)$, $(M_n)$ or $(N_n)$.
Denote the martingale increment $\eps_{n+1} = Y_{n+1}-\gamma_nY_n$. We obtain that
\begin{align*}
    \Delta\cM_{n+1} & =\cM_{n+1}-\cM_n \\
                    & = \begin{pmatrix}
                        S_{n+1}-S_n + \frac{a(\beta+1)}{\beta-a(\beta+1)}\big(\frac{Y_{n+1}}{\mu_{n+1}}-\frac{Y_n}{\mu_n}\big) \\
                            { a_{n+1}Y_{n+1}- a_nY_n}
                    \end{pmatrix}\\
                    & = \begin{pmatrix}\big(1+\frac{a(\beta+1)}{\beta-a(\beta+1)}\big)X_{{n+1}} - \frac{a(\beta+1)}{(\beta-a(\beta+1))\mu_{n+1}}\frac{\beta}{n}Y_n \\
                    { a_{n+1}\eps_{n+1}}
                    \end{pmatrix} \\
                    & = \begin{pmatrix}\frac{\beta}{(\beta-a(\beta+1))\mu_{n+1}}\big(Y_n + X_{{n+1}}\mu_{n+1} - (\gamma_n-1)Y_n\big) \\
                    { a_{n+1}\eps_{n+1}}.
                    \end{pmatrix} 
\end{align*}
Consequently
\begin{equation*}
    \Delta\cM_{n+1}  = \comb{\frac{\beta}{(\beta-a(\beta+1))\mu_{n+1}}}{ a_{n+1}}\eps_{n+1}.
\end{equation*}

We also obtain that
\begin{align}
\label{AERW-ESPC-eps2}
    \E[\eps_{n+1}^2\mid \cF_n]&= \E[Y_{n+1}^2\mid \cF_n]-\gamma_n^2Y_n^2 \nonumber\\
    & = Y_n^2 + 2 (\gamma_n-1)Y_n^2 + \mu_{n+1}^2 - \gamma_n^2 Y_n^2 \nonumber\\
    & = \mu_{n+1}^2 - (\gamma_n-1)^2Y_n^2.
\end{align}
Therefore, we deduce that
\begin{align*}
    \E\big[(\Delta\cM_{n+1})(\Delta\cM_{n+1})^T &\mid \cF_n\big]  = \\
    &(\mu_{n+1}^2 - (\gamma_n-1)^2Y_n^2)\begin{pmatrix}\big(\frac{\beta}{(\beta-a(\beta+1))\mu_{n+1}}\big)^2  & \frac{\beta a_{n+1}}{(\beta-a(\beta+1))\mu_{n+1}} \\ \frac{\beta a_{n+1}}{(\beta-a(\beta+1))\mu_{n+1}} &  a_{n+1}^2 \end{pmatrix}.
\end{align*}
We are now able to compute the quadratic variation of $\cM_n$
\begin{equation}
\label{AERW-QV-MMN}
    \langle \cM\rangle_n = \sum_{k=0}^{n-1}
    \begin{pmatrix} \big(\frac{\beta}{\beta-a(\beta+1)}\big)^2  &  \frac{\beta a_{k+1}\mu_{k+1}}{\beta-a(\beta+1)} \\ \frac{\beta a_{k+1}\mu_{k+1}}{\beta-a(\beta+1)} & ( a_{k+1}\mu_{k+1})^2 \end{pmatrix} - \xi_n
\end{equation}
where 
\begin{equation*}
\xi_n = \sum_{k=0}^{n-1} (\gamma_k-1)^2Y_k^2\begin{pmatrix} \big(\frac{\beta}{(\beta-a(\beta+1))}\big)^2  &  \frac{\beta a_{k+1}\mu_{k+1}}{(\beta-a(\beta+1))} \\ \frac{\beta a_{k+1}\mu_{k+1}}{(\beta-a(\beta+1))} & ( a_{k+1}\mu_{k+1})^2 \end{pmatrix}.
\end{equation*}

Hereafter, we immediatly deduce from \eqref{AERW-QV-MMN} that 
\begin{equation}
\label{AERW-QV-M}
    \cMc_n = \sum_{k=1}^n ( a_k\mu_k)^2 - \zeta_n \quad\text{where}\quad \zeta_n = \sum_{k=1}^n a_k^2(\gamma_k-1)^2Y_k^2
\end{equation}
 and 
\begin{equation}
\label{AERW-QV-N}
    \langle N\rangle_n = \Big(\frac{\beta}{\beta-a(\beta+1)}\Big)^2n.
\end{equation}
The asympotic behavior of $M_n$ is closely related to the one of
\begin{equation}
    \label{AERW-DEF-sumVN}
    w_n = \sum_{k=1}^n ( a_k\mu_k)^2
\end{equation}
as one can observe that we always have $\cMc_n \leq w_n$ and that $\zeta_n$ is negligeable when compared to $w_n$. Consequently, it follows from the definitions of $( a_n)$ and $(\mu_n)$that we have three regimes of behavior for $(M_n)$. In the diffusive regime where is $p<\frac{4\beta+3}{4(\beta+1)}$ or $a<1-\frac{1}{2(\beta+1)}$,
\begin{equation}
\label{AERW-VN-DIF}
    \lim_{n\to\infty} \frac{w_n}{n^{1-2(a(\beta+1)-\beta)}}=\ell \quad\text{where}\quad
    \ell=\frac{1}{1+2(\beta-a(\beta+1))}\Big(\frac{\Gamma(a(\beta+1)+1)}{\Gamma(\beta+1)}\Big)^2.
\end{equation}
In the critical regime where $p=\frac{4\beta+3}{4(\beta+1)}$ or $a=1- \frac{1}{2(\beta+1)}$,
\begin{equation}
\label{AERW-VN-CR}
    \lim_{n\to\infty} \frac{w_n}{\log n}=\Big(\frac{\Gamma(\beta+1+\frac{1}{2})}{\Gamma(\beta+1)}\Big)^2.
\end{equation}
In the superdiffusive regime where $p>\frac{4\beta+3}{4(\beta+1)}$ or $a>1-\frac{1}{2(\beta+1)}$,
\begin{equation}
\label{AERW-VN-SDIF}
    \lim_{n\to\infty} {w_n}= \sum_{k=1}^\infty \Big(\frac{\Gamma(a(\beta+1)+1)\Gamma(k+\beta)}{\Gamma(k+a(\beta+1))\Gamma(\beta+1)}\Big)^2 < +\infty.
\end{equation}

\section{Proofs of the main results}
\label{AERW-S-proofs}

\subsection{The diffusive regime}

%

\begin{lemma}
\label{AERW-lemMV} 
Let $(V_n)$ be the sequence of positive definite diagonal matrices of order $2$ given by
\begin{equation}
\label{AERW-DEF-VN}
         V_n = \frac{1}{\sqrt{n}}\begin{pmatrix} 1 & 0 \\0 & \frac{a(\beta+1)}{\beta-a(\beta+1)}( a_n\mu_n)^{-1}\end{pmatrix}.
\end{equation}
Let $v = \begin{pmatrix}1 \\-1\end{pmatrix}$ such that
\begin{equation}
\label{AERW-SN-vVNMMN}
        v^T V_n\cM_n = \frac{S_n}{\sqrt{n}}.
\end{equation}The quadratric variation of $\langle \cM\rangle_n$ satisfies in the diffusive regime where is $a<1-\frac{1}{2(\beta+1)}$,
\begin{equation}
    \label{AERW-lim-VMM}
        \lim_{n\to\infty} V_n \langle \cM\rangle_n V_n^T = V \quad\text{a.s.}
    \end{equation}
    where the matrix $V$ is given by
\begin{equation}
\label{AERW-DEF-V}
        V = \frac{1}{(\beta-a(\beta+1))^2}\begin{pmatrix}\beta^2 & \frac{a\beta}{1-a}
    \\ \frac{a\beta}{1-a} & \frac{a^2(\beta+1)^2}{1+2\beta-2a(\beta+1)}. \end{pmatrix}
\end{equation}
\end{lemma}

\begin{remark}
    Following the same steps as in the proof of Lemma \ref{AERW-lemMV}, we find that in the critical regime $a = 1-\frac{1}{2(\beta+1)}$, the sequence of normalization matrices $(V_n)$ has to be replaced by 
\begin{equation}
\label{AERW-DEF-WN}
    W_n = \frac{1}{\sqrt{n\log n}}\begin{pmatrix} 1 & 0 \\0 & (2\beta+1)( a_n\mu_n)^{-1}\end{pmatrix}.
\end{equation}
The limit matrix $V$ also need to be replaced by 
\begin{equation}
\label{AERW-DEF-W}
        W = (2\beta+1)^2\begin{pmatrix} 0 & 0 \\ 0 & 1\end{pmatrix}.
\end{equation}
\end{remark}
\begin{proof}{Lemma}{\ref{AERW-lemMV}}
We obtain from Theorem \ref{AERW-THM-CVPS-DIF}, equations \eqref{AERW-lim-alphaN} and \eqref{AERW-VN-DIF} that
\begin{align*}
    \lim_{n\to\infty}& V_n \langle \cM\rangle_n V_n^T
    \\
    & = \lim_{n\to\infty} \frac{1}{n} \begin{pmatrix} \sum_{k=0}^{n-1} \big(\frac{\beta}{(\beta-a(\beta+1))}\big)^2 & \frac{a(\beta+1)\beta}{(\beta-a(\beta+1))^2 a_n\mu_n}\sum_{k=0}^{n-1} a_{k+1}\mu_{k+1} \\ \frac{a(\beta+1)\beta}{(\beta-a(\beta+1))^2 a_n\mu_n}\sum_{k=0}^{n-1} a_{k+1}\mu_{k+1} & \big(\frac{a(\beta+1)}{(\beta-a(\beta+1)) a_n\mu_n}\big)^2\sum_{k=0}^{n-1}( a_{k+1}\mu_{k+1})^2 \end{pmatrix} \\
    & = \frac{1}{(\beta-a(\beta+1))^2}\begin{pmatrix}\beta^2 & \frac{a(\beta+1)\beta}{\beta+1-a(\beta+1)}
    \\ \frac{a(\beta+1)\beta}{\beta+1-a(\beta+1)} & \frac{a(\beta+1)^2}{2(\beta-a(\beta+1))+1} \end{pmatrix}
\end{align*}
which is exactly what we wanted to prove.
\end{proof}

\begin{proof}{Theorem}{\ref{AERW-THM-CVPS-DIF}}
We shall make extensive use of the strong law of large numbers for martingales given, e.g. by theorem 1.3.24 of \cite{Duflo1997}. First, we have for $(M_n)$ that for any $\gamma>0$,
\begin{equation*}
    M_n^2 = O\big((\log w_n)^{1+\gamma}w_n\big)\quad\text{a.s.}
\end{equation*}
which by definition of $M_n$ and as $ a_n$ is asymptotically equivalent to $n^{-a(\beta+1)}$ and $w_n$ is asymptotically equivalent to $n^{1+2(\beta-a(\beta+1))}$ ensures that
\begin{equation*}
    \frac{Y_n^2}{n^2} = O\Big((\log n)^{1+\gamma}\frac{n^{1+2(\beta-a(\beta+1))}}{n^{2(1-a(\beta+1))}}\Big)\quad\text{a.s.}
\end{equation*}
Finally as $\mu_n$ is asymptotically equivalent to $n^{\beta}$, we obtain that
\begin{equation*}
    \frac{Y_n^2}{(\mu_nn)^2} = O\Big(\frac{(\log n)^{1+\gamma}}{n}\Big)\quad\text{a.s.}
\end{equation*}
which reduces to 
\begin{equation}
\label{AERW-lim-YN}
    \lim_{n\to\infty} \frac{Y_n}{\mu_nn}=0 \quad\text{a.s.}
\end{equation}
We now focus our attention on $(N_n)$. By the same token as before, we have that for any $\gamma>0$,
\begin{equation*}
    N_n^2 = O\big((\log n)^{1+\gamma}n\big)\quad\text{a.s.}
\end{equation*}
which by definition of $(N_n)$ gives us
\begin{equation*}
    \frac{\big(S_n-\frac{a(\beta+1)}{\beta-a(\beta+1)}\mu_n^{-1}Y_n\big)^2}{n^2} = O\Big(\frac{(\log n)^{1+\gamma}}{n}\Big)\quad\text{a.s.}
\end{equation*}
and we conclude that
\begin{equation}
\label{AERW-lim-SYN}
    \lim_{n\to\infty} \frac{S_n}{n}-\frac{a(\beta+1)}{\beta-a(\beta+1)}\frac{Y_n}{\mu_nn}=0 \quad\text{a.s.}
\end{equation}
This achieves the proof of Theorem \ref{AERW-THM-CVPS-DIF} as the convergences \eqref{AERW-lim-YN} and \eqref{AERW-lim-SYN} hold almost surely.
\end{proof}

\begin{proof}{Theorem}{\ref{AERW-CLFT-AERW-DR}}
In order to apply Theorem A.2 from \cite{Laulin2021}, we must verify that (H.1), (H.2) and (H.3) are satisfied. \medskip \\
\textbf{(H.1)} We have from \eqref{AERW-lim-VMM} and the fact that $a_{\floor{nt}}$ is asymtotically equivalent to $t^{-a(\beta+1)}a_n$ that
\begin{equation*}
V_n \langle \cM \rangle_{\floor{nt}} V_n^T {\underset{n\to\infty}{\longrightarrow}} V_t \quad\text{a.s.}
\end{equation*}
where
\begin{equation*}
V_t = \dfrac{1}{(\beta-a(\beta+1))^2}\begin{pmatrix}\beta^2 t & \dfrac{a\beta}{1-a} t^{1+\beta-a(\beta+1)}
    \\ \dfrac{a\beta}{1-a} t^{1+\beta-a(\beta+1)} & 
    \dfrac{a^2(\beta+1)^2}{1+2\beta-2a(\beta+1)} t^{1+2\beta-2a(\beta+1)} \end{pmatrix}.
\end{equation*}
\medskip \\
\textbf{(H.2)} In order to verify that Lindeberg's condition is satisfied, we start by deducing from \eqref{AERW-DEF-M-N} together with \eqref{AERW-DEF-MMN} and $V_n$ given by \eqref{AERW-DEF-VN}
that for all $1 \leq k \leq n$
\begin{equation*}
    V_n \Delta\cM_{k}  = \frac{1}{(\beta-a(\beta+1))\sqrt{n}\mu_{n}}\comb{\beta\frac{\mu_n}{\mu_k}}{a\frac{a_k}{a_n}}\eps_{k}
\end{equation*}
which implies that 
\begin{equation}
\label{AERW-VNMMN-BOUND2}
\| V_n \Delta \cM_k \|^2 
= \frac{1}{(\beta-a(\beta+1))^2n}\Big(\frac{\beta^2}{\mu_{k}^2}+\frac{a^2a_k^2}{(a_n\mu_n)^2}\Big)\eps_{k}^2.
\end{equation}
Consequently, we obtain that for all $\eps >0$,
\begin{equation}
\sum_{k=1}^n \E\big[\|V_n \Delta \cM_k \|^2 \1_{\{\|V_n\Delta \cM_k \|>\eps\}}\mid \cF_{k-1}\big]
    \leq  \frac{1}{\eps^2}\sum_{k=1}^n \E\big[\|V_n \Delta \cM_k \|^4\mid \cF_{k-1}\big]. 
\label{AERW-LINDEBERG-DR1} 
\end{equation}
It follows from \eqref{AERW-lim-alphaN} that
\begin{align*}
a_n^{-2}\sum_{k=1}^n a_k^2  = O(n) \quad\text{and}\quad
a_n^{-4}\sum_{k=1}^n a_k^4  = O(n).
\end{align*}
Hence, using that the sequence $(\eps_n)$ is bounded
\begin{equation}
\label{AERW-SUPEPS}
\sup_{1\leq k\leq n} |\eps_k| \leq \sup_{1\leq k\leq n} (\beta+2)\mu_{k} \leq (\beta+2)\mu_{n} \quad\text{a.s.}
\end{equation}
we find that
\begin{equation*}
\sum_{k=1}^n \E\big[\|V_n \Delta \cM_k \|^4 \mid \cF_{k-1}\big] = O\Big(\frac{1}{n}\Big)\quad\text{a.s.}
\end{equation*}
which ensures that Lindeberg's condition (H.2) holds almost surely, that is
for all $\eps >0$,
\begin{equation}
\label{AERW-LINDEBERG-DR4}
\lim_{n \rightarrow \infty} \sum_{k=1}^n \E\big[\|V_n \Delta \cM_k \|^2 \1_{\{\|V_n\Delta \cM_k \|>\eps\}}\mid \cF_{k-1}\big]= 0 \quad \text{a.s.}
\end{equation} 

Since $V_nV_{\floor{nt}}^{-1}$ converges, we immediatly obtain that
\begin{align*}
\lim_{n\to\infty}\sum_{k=1}^{\floor{nt}} \E\bigl[\|V_n \Delta \cM_k \|^2 \1_{\{\|V_n\Delta \cM_k \|>\eps\}}\mid \cF_{k-1}\bigr] & \leq \lim_{n\to\infty}\sum_{k=1}^{\floor{nt}} \E\bigl[\|V_{n} \Delta \cM_k \|^4 \bigr] \\
& \leq\lim_{n\to\infty}\sum_{k=1}^{\floor{nt}} \E\bigl[\|(V_nV_{\floor{nt}}^{-1})V_{\floor{nt}} \Delta \cM_k \|^4\bigr] \\
& = 0\quad\text{a.s.}
\end{align*}

\textbf{(H.3)} In this particular case, we have $V_t = t K_1 + t^{\alpha_2} K_2 + t^{\alpha_3} K_3$ where 
\begin{equation*}
{\alpha_2} = {1-a(\beta+1)}>0 \quad\text{and}\quad {\alpha_3}={1-2a(\beta+1)}>0
\end{equation*} 
as $a<1-\frac{1}{2(\beta+1)}$, and the matrix are symmetric
\begin{equation*}
K_1 =\dfrac{\beta^2}{(\beta-a(\beta+1))^2}\begin{pmatrix}1 & 0 \\ 0 & 0\end{pmatrix}, \
K_2 =\dfrac{a\beta}{(1-a)(\beta-a(\beta+1))^2}\begin{pmatrix}0 & 1 \\ 1 & 0\end{pmatrix},
\end{equation*}
\begin{equation*}
K_3 =\dfrac{a^2(\beta+1)^2}{(1+2\beta-2a(\beta+1))(\beta-a(\beta+1))^2}\begin{pmatrix}0 & 0 \\ 0 & 1\end{pmatrix}.
\end{equation*}
Consequently, we obtain that
\begin{equation*}
\big(V_n \cM_{\floor{nt}}, \ {t \geq 0}\big) \Longrightarrow \big(\cB_{t}, \ {t \geq 0}\big)
\end{equation*}
where $\cB$ is defined as in Theorem A.2 from \cite{Laulin2021}.
Finally,  using the fact that $S_{\floor{nt}}$ is asymptotically equivalent to $N_{\floor{nt}}+t^{\beta-a(\beta+1)}\frac{a(\beta+1)}{\beta-a(\beta+1)}(\mu_n a_{n})^{-1}M_{\floor{nt}}$, and  
 multiplying by $u_t=\comb{1}{t^{a(\beta+1)-\beta}}$, we conclude
\begin{equation}
\label{AERW-CLTF-MERW}
\big(\frac{1}{\sqrt{n}} S_{\floor{nt}}, \ {t \geq 0}\big) \Longrightarrow \big(W_t, \ {t \geq 0}\big)
\end{equation}
where $W_t =u_t^T \cB_t$. It only remains to compute the covariance function of $(W_t)$ that is for $0\leq s \leq t$
\begin{align*}
\E\big[W_sW_t\big] & = u_s^T\E\big[\cB_s \cB_t^T\big]u_t \\
& = u_s^TV_su_t \\
& = u_s^T\big(sK_1+ s^{1+\beta-a(\beta+1)}K_2+s^{1+2\beta-2a(\beta+1)}K_3)u_t \\
& = \frac{\beta^2}{(\beta-a(\beta+1))^2}s +\frac{a\beta s^{1+\beta-a(\beta+1)}}{(1-a)(\beta-a(\beta+1))^2}(s^{a(\beta+1)-\beta}+t^{a(\beta+1)-\beta})\\
&\quad +\frac{a^2(\beta+1)^2}{(1+2\beta-2a(\beta+1))(\beta-a(\beta+1))^2}s^{1+2\beta-2a(\beta+1)}(st)^{a(\beta+1)-\beta} \\
&= 
\frac{a(1+\beta)(1 - a)+a\beta}{(2(\beta+1)(1-a)-1)(a-\beta(1-a))(1-a)}s\Big(\frac{t}{s}\Big)^{a-\beta(1-a)} \nonumber \\
& \quad+ \frac{\beta}{(\beta(1-a)-a)(1-a)}s.
\end{align*}
\end{proof} 

\begin{proof}{Theorem}{\ref{AERW-THM-LFQ-DIF}}
We need to check that all the hypotheses of Theorem A.3 in \cite{Laulin2021} are satisfied. Thanks to Lemma \ref{AERW-lemMV}, hypothesis $(\textnormal{H.1})$ holds almost surely. We also immediately obtain from \eqref{AERW-LINDEBERG-DR4} that $(\textnormal{H.2})$ is verified almost surely when $t=1$.

Hereafter, we need to verify $(\textnormal{H.4})$ is satisfied in the special case $\beta=2$ that is
\begin{equation*}
\label{LFQ-DIF1}
\sum_{n=1}^{\infty} \frac{1}{\bigl(\log  (\det V_{n}^{-1})^2\bigr)^{2}}\E\big[\|V_{n} \Delta \cM_{n}\|^{4}\big|\cF_{n-1}\big]<\infty \quad \text{a.s.}
\end{equation*}
We immediately have from \eqref{AERW-DEF-VN}
\begin{equation}
\label{AERW-DET-VN}
\det V_n^{-1} = \frac{\beta-a(\beta+1)}{a(\beta+1)}a_n\mu_n\sqrt{n}.
\end{equation} 
Hence, we obtain from \eqref{AERW-lim-alphaN} and \eqref{AERW-DET-VN} that
\begin{equation}
\label{LOG-DET-VN}
\lim_{n \rightarrow \infty}  \frac{\log  (\det V_{n}^{-1})^2 }{\log n} 
= 1+2\beta-2a(\beta+1).
\end{equation}
Therefore, we can replace $\log  (\det V_{n}^{-1})^2$ by $\log n$ in \eqref{LFQ-DIF1}. Hereafter, we obtain from \eqref{AERW-VNMMN-BOUND2} and \eqref{AERW-SUPEPS} that
\begin{align}
\label{LFQ-DIF2}
\sum_{n=2}^{\infty} \frac{1}{(\log n)^2}\E\big[\|V_{n} \Delta \cM_{n}\|^{4}\big|\cF_{n-1}\big]
& = O \Big( \sum_{n=1}^{\infty} \frac{1}{(n\log n)^2}  \Big). 
\end{align}
Thus, \eqref{LFQ-DIF2} guarentees that $(\textnormal{H.4})$ is verified.
We are now going to apply the quadratic strong law given by Theorem A.3 in \cite{Laulin2021}. We get from equation \eqref{LOG-DET-VN} that
\begin{equation}
\label{LFQ-DIF3}
\lim_{n \rightarrow \infty}  \frac{1}{\log n}\sum_{k=1}^n
\Big(\frac{(\det V_{k})^2 - (\det V_{k+1})^2}{(\det V_k)^2}\Big)V_k\cM_k\cM_k^T V_k^T =
\big(1+2\beta-2a(\beta+1)\big) V \quad \text{a.s.}
\end{equation}
However, we obtain from \eqref{AERW-lim-alphaN} and \eqref{AERW-DET-VN} that
\begin{equation}
\label{lim-DET-VN}
\lim_{n \rightarrow \infty} n\Big(\frac{(\det V_{n})^2 - (\det V_{n+1})^2}{(\det V_n)^2}\Big)=1+2\beta-2a(\beta+1).
\end{equation}
Finally, we can deduce from \eqref{AERW-SN-vVNMMN}, \eqref{LFQ-DIF3} and \eqref{lim-DET-VN} that
\begin{equation}
\label{LFQ-DIF4}
\lim_{n \rightarrow \infty}  \frac{1}{\log n}\sum_{k=1}^n\frac{S_k^2}{k^2} 
= v^TVv \quad \text{a.s.}
\end{equation}
which, together with
\begin{equation}
\label{COV-vVv}
v^TVv = \frac{2\beta+1-a}{(1-a)(1+2\beta-2a(\beta+1))}
\end{equation}
completes the proof of Theorem \ref{AERW-THM-LFQ-DIF}.
\end{proof}

\subsection{The critical regime}
The proofs of Theorems {\ref{AERW-THM-CVPS-CR}
and \ref{AERW-THM-CLFT-CR} follows essentially the same lines as the ones in the diffusive regimes, provided one exchange $V_n$ with $W_n$. Hence, they shall not be explicited here.
\begin{proof}{Theorem}{\ref{AERW-THM-LFQ-LIL-CR}}
The proof of the quadratic strong law \eqref{AERW-QUAD-CR} is left to the reader as it 
follows essentially the same lines as that of \eqref{AERW-QUAD-DIF}. The only minor change is that the matrix $V_n$ has to be replaced by the matrix $W_n$ defined in \eqref{AERW-DEF-WN}.
We shall now proceed to the proof of the law of iterated logarithm
given by \eqref{AERW-LIL-CR}. On the one hand, it follows from \eqref{AERW-lim-alphaN} and \eqref{AERW-VN-DIF} that 
\begin{equation}
\label{CONDLIL-CR}
\sum_{n=1}^{+\infty} \frac{a_n^4}{w_n^2} < \infty.
\end{equation}
Moreover, we have from \eqref{AERW-QV-M} and \eqref{AERW-QV-N} that 
\begin{equation*}
\lim_{n\to\infty}\frac{\langle M\rangle_n}{w_n} =  1\quad\text{a.s.}
\quad\text{and}\quad
\lim_{n\to\infty}\frac{\langle N\rangle_n}{n} = \Big(\frac{\beta}{\beta-a(\beta+1)}\Big)^2 \quad\text{a.s.}
\end{equation*}
Consequently, we deduce from the law of iterated logarithm for martingales due to Stout \cite{Stout1974}, see also Corollary 6.4.25 in \cite{Duflo1997},  
that $(M_n)$ satisfies when $a=1-1/2(\beta+1)$
\begin{align*}
\limsup_{n \rightarrow \infty} \frac{M_n}{(2 w_n \log \log w_n)^{1/2}}  & = -\liminf_{n \rightarrow \infty} \frac{M_n}{(2 w_n \log \log w_n)^{1/2}}  \\
& = 1  \quad\text{a.s.}
\end{align*}
However, as $a_n w_n^{-1/2}$ is asymptotically equivalent to $(n^{2\beta+1}\log n)^{-1/2}$, we immediately obtain from \eqref{AERW-VN-CR} that
\begin{align}
\label{LIL-CR1}
\limsup_{n \rightarrow \infty} \frac{Y_n}{(2 n^{2\beta+1} \log n \log \log \log n)^{1/2}}  & = -\liminf_{n \rightarrow \infty} \frac{Y_n}{(2 n^{2\beta+1} \log n \log \log \log n)^{1/2}} \nonumber  \\
\limsup_{n \rightarrow \infty} \frac{n^{-\beta}Y_n}{(2 n \log n \log \log \log n)^{1/2}}  & = -\liminf_{n \rightarrow \infty} \frac{n^{-\beta}Y_n}{(2 n \log n \log \log \log n)^{1/2}} \nonumber  \\
& = 1  \quad\text{a.s.}
\end{align}
The law of iterated logarithm for martingales also allow us to find that $(N_n)$ satisfies
\begin{align*}
\label{LIL-CR2}
\limsup_{n \rightarrow \infty} \frac{N_n}{(2 n \log \log n)^{1/2}}  & = -\liminf_{n \rightarrow \infty} \frac{N_n}{(2 n \log \log n)^{1/2}}  \\
& = \sqrt{4\beta^2}  \quad\text{a.s.}
\end{align*}
which ensures that
\begin{equation*}
\limsup_{n \rightarrow \infty}  \frac{N_n}{(2 n \log n \log \log \log n)^{1/2}} 
= 0 \quad\text{a.s.} 
\end{equation*} 
Hence, we deduce from \eqref{AERW-SN-MG} and \eqref{LIL-CR1} that
\begin{align*}
\label{LIL-CR3}
\limsup_{n \rightarrow \infty} \frac{S_n}{(2 n\log n \log \log \log n)^{1/2}}  
& = \limsup_{n \rightarrow \infty} \frac{N_n + (2\beta+1)(\mu_n a_n)^{-1}M_n}{(2 n\log n \log \log \log n)^{1/2}} \nonumber \\
& = \limsup_{n \rightarrow \infty} (2\beta+1)\frac{Y_n}{(2 n^{2\beta+1}\log n \log \log \log n)^{1/2}} \nonumber \\
& = -\liminf_{n \rightarrow \infty} (2\beta+1)\frac{Y_n}{(2 n^{2\beta+1}\log n \log \log \log n)^{1/2}} \nonumber  \\
& = -\liminf_{n \rightarrow \infty} \frac{S_n}{(2 n\log n \log \log \log n)^{1/2}}. \nonumber
\end{align*}
Hence, we obtain that 
\begin{align*}
\limsup_{n \rightarrow \infty} \frac{S_n^2}{2 n\log n \log \log \log n}
&= \limsup_{n \rightarrow \infty} (2\beta+1)^2\frac{Y_n^2}{2 n\log n \log \log \log n} \\
&=(2\beta+1)^2
\end{align*}
which immediately leads to \eqref{AERW-LIL-CR}, thus completing the proof of Theorem \ref{AERW-THM-LFQ-LIL-CR}.
\end{proof}

\subsection{The superdiffusive regime}
%
\begin{proof}{Theorem}{\ref{AERW-THM-CVPS-SDIF}}
Hereafter, we shall again make extensive use of the strong law of large numbers for martingales given, e.g. by theorem 1.3.24 of \cite{Duflo1997} in order to prove \eqref{AERW-CVPS-SDIF}. When $a>1-\frac{1}{2(\beta+1)}$, we have from \eqref{AERW-VN-SDIF} that $w_n$ converges. Hence, as $\langle M\rangle_n \leq w_n$, we clealy have that $\langle M\rangle_\infty <\infty$ almost surely and we can conclude that
\begin{equation*}
\lim_{n\to\infty} M_n = M \quad\text{a.s.} \quad\text{where}\quad M = \sum_{k=1}^{\infty} a_k \eps_k
\end{equation*}
which by definition of $M_n$, and as $a_n$ is asymptotically equivalent to $\Gamma(a(\beta+1)+1)n^{-a(\beta+1)}$, ensures that
\begin{equation}
\label{lim-YN-SDIF}
\lim_{n\to\infty}\frac{Y_n}{n^{a(\beta+1)}} = Y \quad\text{a.s.}\quad\text{where}\quad Y = \frac{1}{\Gamma(a(\beta+1)+1)}M.
\end{equation}
Moreover, we still have that for any $\gamma>0$,
\begin{equation*}
N_n^2 = O\big((\log n)^{1+\gamma}n\big)\quad\text{a.s.}
\end{equation*}
which by definition of $N_n$ gives us for all $t\geq 0$
\begin{equation*}
\frac{\big(S_n + \frac{a(\beta+1)}{\beta-a(\beta+1)}(\mu_n)^{-1}Y_n
)^2}{n^{2a(\beta+1)-2\beta}} = O\Big(\frac{(\log n)^{1+\gamma}}{n^{2a(\beta+1)-2\beta-1}}\Big)\quad\text{a.s.}
\end{equation*}
We know that $a>1-\frac{1}{2(\beta+1)}$ in the superdiffusive regime, which ensures that $2a(\beta+1)-2\beta-1>0$. Then, we obtain thanks to \eqref{AERW-lim-muN} and \eqref{AERW-lim-YN} that for all $t\geq 0$
\begin{equation}
\label{lim-SYN-SDIF1}
\lim_{n\to\infty} \frac{S_{\floor{nt}}}{{\floor{nt}}^{a(\beta+1)-\beta}}+\frac{a(\beta+1)}{\beta-a(\beta+1)}\frac{Y_{\floor{nt}}}{{\floor{nt}}^{a(\beta+1)}}=0 \quad\text{a.s.}
\end{equation}
The convergences \eqref{lim-YN-SDIF} and \eqref{lim-SYN-SDIF1} hold almost surely and $\floor{nt}$ is asymptotically equivalent to $nt$ which implies
\begin{equation}
\label{lim-SYN-SDIF-2}
\lim_{n\to\infty} \frac{S_{\floor{nt}}}{n{^a(\beta+1)}}= t^{a(\beta+1)}L_\beta \quad\text{a.s.}
\end{equation} 
Finally, the fact that \eqref{lim-SYN-SDIF-2} holds almost surely ensures that it also holds for the finite-dimensional distributions, and we obtain \eqref{AERW-CVPS-cont-SDIF} with $\Lambda_t = t^{a(\beta+1)}L_\beta$ and $L_\beta =\frac{a(\beta+1)}{a(\beta+1)-\beta}Y$. \medskip

We shall now proceed to the proof of the mean square convergence \eqref{AERW-CVL2-SDIF}. On the one hand, as $M_0=0$ we have from \eqref{AERW-QV-M} that
\begin{equation*}
\E\big[M_n^2\big]=\E\big[\langle M\rangle_n\big]\leq w_n.
\end{equation*}
Hence, we obtain from \eqref{AERW-VN-SDIF} that
\begin{equation*}
\sup_{n\geq 1} \E\big[M_n^2\big] < \infty
\end{equation*}
which ensures that the martingale $(M_n)$ is bounded in $\IL^2$. Therefore, we have the mean square convergence
\begin{equation*}
\lim_{n\to\infty} \E\big[\big|M_n-M\big|^2\big]=0
\end{equation*}
which implies that
\begin{equation}
\label{L2-Y}
\lim_{n\to\infty} \E\Big[\Big|\frac{Y_n}{n^{a(\beta+1)}}-Y\Big|^2\Big]=0.
\end{equation}
On the other hand, for any $n\geq0$, the martingale $(N_n)$ satisfies
\begin{equation*}
\E\big[N_n^2\big]=\E\big[\langle N\rangle_n\big]\leq \Big(\frac{\beta}{\beta-a(\beta+1)}\Big)^2n
\end{equation*}
and since $a(\beta+1)-\beta>\frac{1}{2}$ we obtain
\begin{equation}
\label{L2-N}
\lim_{n\to\infty} \E\Big[\Big|\frac{N_n}{n^{a(\beta+1)-\beta}}\Big|^2\Big]=0.
\end{equation}
Finally, we obtain the mean square convergence \eqref{AERW-CVL2-SDIF} from \eqref{L2-Y} and \eqref{L2-N} and we achieve the proof of Theorem \ref{AERW-THM-CVPS-SDIF}.
\end{proof}

%
\vspace{10ex}
\bibliographystyle{acm}
{\small\setstretch{0.1}\bibliography{bibliography}}

\end{document}